\documentclass[10pt]{amsart}
\usepackage{amsmath}
\usepackage{amssymb}
\usepackage{mathrsfs}		
\usepackage[all]{xy}
\usepackage[colorlinks=true,urlcolor=blue,bookmarks=true,bookmarksopen=true,citecolor=blue,hypertex]{hyperref}
\def \qed{\hfill $\square$ \vspace{0.03in}}

\numberwithin{equation}{section}

\newtheorem{theorem}{Theorem}[section]
\newtheorem{defn}[theorem]{Definition}
\newtheorem{theoremdefn}[theorem]{Theorem(Definition)}
\newtheorem{assumption}[theorem]{Assumption}

\newtheorem{corollary}[theorem]{Corollary}

\newtheorem{lemma}[theorem]{Lemma}
\newtheorem{prop}[theorem]{Proposition}
\newtheorem{remark}[theorem]{Remark}

\def \begineq{\begin{equation}}
\def \endeq{\end{equation}}

\def \n{\noindent}

\def \bb{\mathbb}

\def \mc{\mathcal}
\def \mf{\mathfrak}
\def \ms{\mathscr}

\def \CC{{\bb C}}

\def \JJ{{\bb J}}

\def \RR{{\bb R}}
\def \TT{{\bb T}}
\def \VV{{\bb V}}

\def \ZZ{{\bb Z}}

\def \DDC{{\mc D}}
\def \JJC{{\mc J}}
\def \LLC{{\mc L}}

\def \TTC{{\mc T}}

\def \GGS{{\ms G}}

\def \XXS{{\ms X}}

\def \({\left(}
\def \){\right)}
\def \<{\langle}
\def \>{\rangle}
\def \bar{\overline}

\def \dsum{\oplus}
\def \eps{\epsilon}

\def \inter{\cap}
\def \into{\hookrightarrow}

\def \xto{\xrightarrow}

\def \Ann{{\rm Ann}}
\def \Diff{{\rm Diff}}
\def \Ham{{\rm Ham}}
\def \Hom{{\rm Hom}}
\def \img{{\rm img}}
\def \Jac{{\rm Jac}}
\def \Nij{{\rm Nij}}
\def \Span{{\rm Span}}

\begin{document}

\title{Hamiltonian symmetries and reduction in generalized geometry}
\author{Shengda Hu}
\email{shengda@dms.umontreal.ca}
\address{D\'epartement de Math\'ematiques et de Statistique, 
Universit\'e de  Montr\'eal, CP 6128 succ Centre-Ville, Montr\'eal, QC H3C 3J7, Canada}

\abstract
A closed $3$-form $H \in \Omega^3_0(M)$ defines an extension of $\Gamma(TM)$ by $\Omega^2_0(M)$. This fact leads to the definition of the group of $H$-twisted Hamiltonian symmetries $\Ham(M, \JJ; H)$ as well as Hamiltonian action of Lie group and moment map in the category of (twisted) generalized complex manifold. The Hamiltonian reduction in the category of generalized complex geometry is then constructed. 
The definitions and constructions are natural extensions of the corresponding ones in the symplectic geometry. We describe cutting in generalized complex geometry to show that it's a general phenomenon in generalized geometry that topology change is often accompanied by twisting (class) change. 
\endabstract
\maketitle


\section{Introduction}\label{intro}

Generalized complex structure was introduced by Hitchin \cite{Hitchin} and developed by Gualtieri in his thesis \cite{Gualtieri}. Much more work has been done since. This new structure specializes to complex and symplectic structures on two extremes and it is regarded as the natural category in which to consider constructions pertaining to both, e.g. mirror symmetry.

The group of Hamiltonian diffeomorphisms occupies a prominent place in symplectic geometry. Reduction via group actions and related constructions have been proved fruitful whenever they exist. The existence of such notions and constructions in generalized geometry would definitely be desired (e.g. \cite{Kapustin}). The notion of infinitesimal Hamiltonian symmetry in (untwisted) generalized geometry is given by \cite{Gualtieri}, in remark $5.2$ and proposition $5.3$. 
We obtain the following analogues of these notions from symplectic geometry in (twisted) generalized complex geometry:

\vspace{0.1in}
\noindent
{\bf Theorem \ref{ham:gp}.} {\it Let $(M, \JJ)$ be an $H$-twisted generalized complex manifold, $\GGS_\JJ$ be the subgroup of $H$-twisted generalized symmetries preserving $\JJ$. Then $\Ham(M, \JJ)$, the set of time-$1$ generalized complex symmetries generated by time-dependent Hamiltonian functions is a subgroup of $\GGS_\JJ$.}

\vspace{0.1in}
\noindent
{\bf Theorem \ref{quot:algebroid}.} {\it Let $G$ be a compact Lie group. Suppose $(M, \JJ)$ is an $H$-twisted generalized complex manifold with a Hamiltonian $G$-action with moment map $\mu : M \to \mf g^*$, so that $0$ is a regular value of $\mu$ and the geometrical action of $G$ is free on $M_0 = \mu^{-1}(0)$. Then there is a natural extended complex structure on the reduced space $Q = M_0/ G$.}

\vspace{0.1in}
We describe the content of the article in the following. 
The generalized geometry studies structures on the generalized tangent bundle $\TT M = TM \dsum T^*M$ with the structure of Courant algebroid given by a closed three form $H \in \Omega^3_0(M)$. When $H = 0$, it's shown in \cite{Gualtieri} that the group of symmetries of the Courant algebroid $\TT M$ is given by the semi-direct product
$$\GGS_0 = \Diff(M) \ltimes \Omega^2_0(M), \text{ whose Lie algebra is } \XXS_0 = \Gamma(TM) \dsum \Omega_0^2(M).$$
The Lie algebra structure on $\XXS_0$ is 
\begin{equation}\label{intro:bracket}
[(X, A), (Y, B)] = ([X, Y], \LLC_X B - \LLC_Y A),
\end{equation}
which can be seen as the trivial extension of $\Gamma(TM)$ by the module $\Omega^2_0(M)$:
$$X \circ A = \LLC_X A \text{ for } X \in \Gamma(TM) \text{ and } A \in \Omega^2_0(M).$$
In the presence of twisting $H$, let $\GGS_H$ and $\XXS_H$ be the respective symmetry group and Lie algebra. We show in \S\ref{subg} that $\XXS_H$ is the extension of $\Gamma(TM)$ by the module $\Omega^2_0(M)$ defined by the $2$-cocyle $\alpha_H$
$$\alpha_H(X, Y) = d\iota_Y \iota_X H, \text{ for } X, Y \in \Gamma(TM).$$
Let $\XXS = \Gamma(TM) \dsum \Omega^2(M)$ with the Lie bracket \eqref{intro:bracket}. Then the inclusion 
$$\psi_H : \XXS_H \into \XXS : (X, A) \to (X, A - \iota_X H)$$
is an inclusion of Lie algebra (proposition \ref{ham:transl}).

The inclusion $\psi_H$ comes into play, for example, when we integrate a generalized field, i.e. a section of $\TT M$, into a symmetry of the $H$-twisted Courant algebroid structure:
$$\mf X = X + \xi \in \Gamma(\TT M) \Longrightarrow e^{t\mf X} := e^{t \psi_H(X, d\xi)} := (\lambda_t = e^{tX}, \int_0^t \lambda_s^* (d\xi - \iota_X H) ds).$$
Similar to symplectic geometry, we define (following \cite{Gualtieri}) the generalized Hamiltonian field generated by $f: M \to \RR$ on an ($H$-twisted) generalized complex manifold $(M, \JJ)$ to be
$$\mf X_f = \JJ(df).$$
The generalized Hamiltonian symmetry generated by $f$ is then $e^{t\mf X_f}$ as above. This enables us to define in \S\ref{ham} the notion of generalized Hamiltonian action by a Lie group $G$ and the group of generalized Hamiltonian symmetries $\Ham(M, \JJ)$, completely parallel to the corresponding definitions in symplectic geometry.

We then consider reduction in our Hamiltonian framework in \S\ref{quot}. 
The induced Courant algebroid structure on the reduced manifold $Q$ is again exact, while there is no canonical identification to $\TT Q$ with twisted Courant bracket, unless the group action factors through $\Diff(M)$, in which case, explicit identification can be written down (upto a choice of connection form) (cf. corollary \ref{quot:descend}). One notable fact of the construction is that the twisting form upstairs is \emph{not} required to descend to the reduced manifold. Furthermore, the \v Severa class of the reduced structure (which is the class $[H] \in H^3(M)$ when a splitting is chosen) does not depend on either the choice of the connection form or invariant $B$-field action on the original manifold. 
We also prove some minor facts, such as the reduction of a generalized Calabi-Yau manifold (as in \cite{Hitchin}) is again generalized Calabi-Yau manifold in the same sense; and 
when the group is torus, the twisting form of the reduced structure satisfies a Duistermaat-Heckman type formula in a component of regular values of the moment map. In \S\ref{quot:exple}, we compute the example of $\CC^2\setminus\{(0,0)\}$ with nontrivially twisted generalized complex structure.

With the generalized complex reduction in hand, other related constructions and phenomena exhibit themselves, such as coupling structure, cutting, wall crossing (at least for $S^1$-action) etc. As in the case of symplectic geometry, we may also weaken the condition of free action and instead have orbifold as reduced space. These we postpone to a later work. Here we only describe cutting (\S\ref{cut}), along the line of cutting constructions of \cite{Lerman1, Lerman2, Burns}, to show that operations in the generalized geometry share a common feature, i.e. change of topology is accompanied by change of twisting (class).

In the appendix, we collect various definitions and theorems on Lie and Courant algebroids and some facts on Lie algebra extensions.

The early version of this work appeared on the arXiv along with the works of Lin-Tolman \cite{Lin}, Bursztyn-Cavalcanti-Gualtieri \cite{Bursztyn} as well as Stienon-Xu \cite{Stienon} around the same time.
In \cite{Hitchin2} Hitchin described an example of generalized K\"ahler reduction. Another version of generalized complex reduction is given in \cite{Crainic}, which is based on a notion of generalized holomorphic map. In particular, the fibers over regular values of such maps carry natural complex structure. Another point of view is provided by \cite{Vaisman}, which discusses three related stuctures on $\TT M$.

The reduction described in this article also fits into bigger picture of reduction theory. We give a very brief and extremely incomplete recount of related works in the following. More may be found in the references of the works mentioned below. First of all, our construction is a direct generalization of the Marsden-Weinstein reduction \cite{Marsden1} in symplectic geometry to generalized complex category. As shown in \cite{Gualtieri}, generalized complex structures can also be defined as complex Dirac structures with real index $0$. In fact, any generalized (as well as extended) complex structure provides the manifold with a Poisson structure (see \S\ref{ham:infinitesimal} and the references therein) and the reduction constructed here coincides with the Poisson reduction defined by \cite{Marsden2}, when seen in the Poisson category. More generally, the reduction of Dirac structures (without twisting) was done in the control community, e.g. \cite{Blankenstein1, Blankenstein2}, where integrability is not required while the reduced structure will be integrable if the original one is so. These are the so called Hamiltonian point of view. Another point of view on reduction of Dirac structures comes from the relation with variational principles, see for example \cite{MarsdenY1, MarsdenY2}. The reduction in singular cases for the Poisson and Dirac structures are also known, e.g. \cite{Ortega}, \cite{Blankenstein3}, as well as the excellent book \cite{Ortega1}, where many more references may be found. Group valued moment maps are discussed in \cite{McDuff1} for $S^1$, \cite{Ginzburg} for torus and \cite{Alekseev} for general case.
The related reduction and further theory can be found in for example \cite{Alekseev1, Alekseev2}, \cite{Bott}. 
Reduction of symplectic structure by action of Poisson Lie group was also discussed, e.g. in \cite{Lu}. In there the moment map would have target space the dual Poisson Lie group instead of dual of the Lie algebra. 
Relation among different sorts of reductions of symplectic structures with or without moment map is discussed in the paper \cite{Ortega1} and the references therein. For the relation of symplectic reduction to algebraic geometry, the survey article \cite{Kirwan} and the references therein are excellent sources.
Of other geometrical structures, to name just a few, such as K\"ahler, hyper-K\"ahler \cite{Hitchin4, Lindstrom, Lindstrom1}, Sasakian, locally conformal symplectic or K\"ahler geometries and contact geometry, etc, the various reductions were considered as well. Because of the limited scope of this paper, we only mention some recent developements of these, from which more references can be found: K\"ahler \cite{Huebschmann2, Huebschmann, Huebschmann0, Huebschmann1, Apostolov}, hyper-K\"ahler \cite{Bielawski, Proudfoot, Ornea}, contact \cite{Albert, Geiges, Willett}, Sasakian \cite{Grantcharov, Dragulete}, Vaisman structure \cite{Gini1}, locally conformal K\"ahler \cite{Gini}, locally conformal symplectic \cite{Noda}, complex Poisson \cite{Nunes} etc. The very incomplete list above could not and was not meant to capture the vast literature on and the span of the reduction theory. Instead, it only shows partly how much more the author needs to learn in this fascinating field of mathematics.

\vspace{0.1in}\n
{\bf Acknowledgement.}
I'd like to thank Fran\c{c}ois Lalonde and DMS in Universit\'e de Montr\'eal for generous support and excellent working conditions, which made this work possible. I'd like to thank Tudor Ratiu for his patience and very helpful comments that brought this work into perspective with respect to general reduction constructions. Of course, the omissions and mistakes on the literature are all due to the author's ignorance in the subject. 
I thank Mainak Poddar and B. Doug Park for invitation to Waterloo, when many enlightening talks were given on generalized complex geometry during the workshop of mirror symmetry at the Perimeter Institute, especially those of Hitchin and Gualtieri, and for helpful discussions that kick-started this project. I thank Sam Lisi for helpful discussions. 
I'd like to thank Yi Lin for initial discussions and providing with an early draft of their paper \cite{Lin}. The current work is partially inspired by their paper. 
I would also like to thank the referee for many valuable suggestions. 
Special thanks go to my family: Aihong, Henry and Catherine for their support and understanding.

\vspace{0.1in}\n
{\bf A note on the notations:}
There are lots of brackets in the following. We did not try to make them all look different, for which there would be too much clutter of notations. Instead, except in the appendix where various bundles are involved and we distinguish them by bundle subscripts, we only make the distinction of $H$-twisted brackets by adding subscript $H$. As to the spaces on which the brackets are defined, it should be clear from the context.

\section{Generalized symmetries}\label{subg}

\subsection{Symmetry group of Courant brackets}
\label{subg:gendef} 
The symmetry group of $(\TT M, \<,\>)$ we'll consider is the \emph{group of generalized symmetries} $\GGS = \Diff(M) \ltimes \Omega^2(M)$, whose action on $\TT M$ is defined as \emph{push-forward} by the following, where $X \in TM$ and $\xi \in T^*M$:
\begineq\label{subg:sym}(\lambda, \alpha) \circ (X + \xi) = \lambda_* X + (\lambda^{-1})^*(\xi + \iota_X \alpha), \text{ where } \lambda \in \Diff(M), \alpha \in \Omega^2(M).
\endeq
We recall the definition of $H$-twisted Courant bracket on $\Gamma(\TT M)$:
$$[X+\xi, Y+\eta]_H = [X, Y] + \LLC_X \eta - \LLC_Y \xi - \frac{1}{2}d(\iota_X\eta - \iota_Y\xi) + \iota_Y \iota_X H.$$
The action of $(\lambda, \alpha)$ on the twisted Courant bracket gives the following:
\begin{equation}\label{subg:bracact}
(\lambda, \alpha) \circ [X+\xi, Y+\eta]_H = [(\lambda, \alpha)\circ (X+\xi), (\lambda, \alpha)\circ (Y + \eta)]_{(\lambda^{-1})^*(H-d\alpha)}.
\end{equation}
The $2$-form $\alpha$ is the \emph{$B$-field} of the generalized symmetry $(\lambda, \alpha)$. The action of $(id, \alpha)$ is also called a \emph{$B$-field transformation}.
Another convention for the action of $\GGS$ is also valid, i.e. $(\lambda, \alpha) \circ (X + \xi) = \lambda_* X + (\lambda^{-1})^*\xi + \iota_{\lambda_*X} \alpha$, while they give the same infinitesimal action. 
We write down the composition law (in our chosen convention for the action):
\begineq\label{subg:comp}
(\lambda, \alpha)\cdot(\mu, \beta) = (\lambda\mu, \mu^*\alpha + \beta).
\endeq
We have the extension sequence of groups:
$$0 \to \Omega^2(M) \to \GGS \xto{\pi_1} \Diff(M) \to 1,$$

Let $\{(\lambda_t, \alpha_t)\}$ be an $1$-parameter subgroup of $\GGS$, then $\{\lambda_t\}$ is an $1$-parameter subgroup of diffeomorphisms generated by a vector field $X$. Let $A = \left.\frac{d\alpha_t}{dt}\right|_{t = 0} \in \Omega^2(M)$ 
then we have
\begin{equation}\label{subg:exp}
(\lambda_t, \alpha_t) = e^{t(X, A)} := (e^{tX}, \int_0^t \lambda_s^* A ds).
\end{equation}
The Lie algebra $\XXS$ of $\GGS$ is then $\Gamma(TM) \dsum \Omega^2(M)$ with the bracket:
\begineq\label{subg:bracket}
[(X, A), (Y, B)] = ([X, Y], \LLC_X B - \LLC_Y A).
\endeq
In other words, $(\XXS, [,])$ is the trivial extension of the standard Lie algebra $\Gamma(TM)$ by $\Omega^2(M)$, where the module structure is given by $X \circ \omega = \LLC_X \omega$ (cf. appendix \S\ref{app:coho}):
$$ 0 \to \Omega^2(M) \to \XXS \to \Gamma(TM) \to 0.$$
More generally, 
let $\{\tilde \lambda_t = (\lambda_t, \alpha_t)\}$ be a smooth path, starting from the identity, in $\GGS$ and $\{X_t\}$ the time-dependent vector fields generating $\{\lambda_t\}$. Then the \emph{infinitesimal symmetry} generating $\{\tilde \lambda_t\}$ is the path $\{(X_t, A_t)\}$ in $\XXS$, where:
\begin{equation}\label{ham:time} 
A_t = (\lambda_t^{-1})^*\dot \alpha_t, \text{ or equivalently } \alpha_t = \int_0^t \lambda_s^*A_s ds.
\end{equation}
The infinitesimal action of $(X, A) \in \XXS$ on $\TT M$ is given by 
\begineq\label{subg:vect}
(X, A) \circ (Y+\eta) = -[X, Y] - \LLC_X \eta + \iota_Y A.
\endeq

A $(k+1)$-form $\rho \in \Omega^{k+1}(M)$ defines a $\Omega^2_0(M)$-valued $k$-cochain $\alpha_\rho$ on the Lie algebra $\Gamma(TM)$:
$$\alpha_\rho(X_1\wedge \ldots \wedge X_k) = d\iota_{X_k}\ldots\iota_{X_1} \rho.$$
We have $\partial_k \alpha_\rho = \alpha_{d\rho}$ where $\partial_k$ is the differential of the Lie algebra cohomology.
Thus, we obtain the following map 
$$\alpha_\bullet: H^{k+1}(M) \to H^k(\Gamma(TM), \Omega^2_0(M)) : [\rho] \mapsto [\alpha_\rho].$$
\begin{defn}\label{ham:twistbrac}
Given $H \in \Omega^3_0(M)$, the Lie algebra $\XXS_H$ is 
the extension of $\Gamma(TM)$ by $\Omega^2_0(M)$ defined by the $2$-cocycle $\alpha_H$. Equivalently, $\XXS_H$ is $\XXS_0 := \Gamma(TM) \dsum \Omega^2_0(M)$ with the following \emph{$H$-twisted Lie bracket} (cf. \eqref{app:lieextbrac}):
\begineq\label{ham:trans}
[(X, A), (Y, B)]_H = ([X, Y], \LLC_X B - \LLC_Y A + d\iota_Y\iota_X H).
\endeq
\end{defn}
\begin{remark}\rm{$\alpha_H$ can also be viewed as a $\Omega^2(M)$-valued $2$-cochain, which defines an extension of Lie algebra $\Gamma(TM)$ by $\Omega^2(M)$. This extension is obviously trivial since $H$ defines another $\Omega^2(M)$-valued $1$-cochain $\gamma_H$ on $\Gamma(TM)$ by $\gamma_H(X) = -\iota_X H$ and we can see that $\partial_1 \gamma_H = \alpha_H$.
}\end{remark}

Let $\GGS_H \subset \GGS$ be the subgroup of generalized symmetries that preserve the Courant bracket $[,]_H$ and $(\lambda_t, \alpha_t)$ a path in $\GGS_H$. Let $(X_t, A_t)$ be the infinitesimal symmetry, then \eqref{subg:bracact} gives
\begin{equation}\label{ham:preserve}
\LLC_{X_t} H + dA_t = d(\iota_{X_t} H + A_t) = 0.
\end{equation}
\begin{prop}\label{ham:transl}
Consider the linear map $\psi_H : \XXS \to \XXS : (X, A) \mapsto (X, A - \iota_X H)$. Then
$$[\psi_H(X, A), \psi_H(Y, B)]_{H'-H} = \psi_H[(X, A), (Y, B)]_{H'}$$
and $\psi_H : (\XXS_0, [,]_H) \to (\XXS, [,])$ is an inclusion of Lie algebra.
\end{prop}
{\it Proof:}
Straightforward computation shows the equality, from which the last statement on $\psi_H$ follows.
\qed

For $(X, A) \in \XXS_0$, the \emph{$H$-twisted infinitesimal symmetry} generated by $(X,A)$ is $\psi_H(X,A) = (X, A - \iota_X H)$.
From \eqref{ham:preserve}, we see that the $1$-parameter subgroup $\{e^{t(X, A- \iota_X H)}\}$ lies in $\GGS_H$. Instead of \eqref{subg:vect}, the $H$-twisted infinitesimal action is
\begineq\label{subg:twistinf}
(X, A) \circ_H (Y+ \eta) = -[X, Y] - \LLC_X \eta + \iota_Y A - \iota_Y\iota_X H.
\endeq
From now on, we will also use $\XXS_H$ to denote the image of $\XXS_H$ under the embedding $\psi_H$.

\subsection{Generalized complex structures} The data $(\<,\>, [,]_H)$ on $\TT M$ together with the natural projection $a$ to $TM$ defines a structure of Courant algebroid as in definition \ref{app:courant}. We then rephrase the definition of generalized complex structure given in \cite{Hitchin, Gualtieri} as follows:
\begin{defn}\label{subg:extend}
The Courant algebroid $(\TT M, \<,\>, [,]_H, a)$ will be called an \emph{extended tangent bundle} and denoted $\TTC M$ when we forget the paticular splitting into the direct sum $TM \dsum T^*M$. The bracket $[,]_H$ will then be denoted simply as $[,]$.
An \emph{extended almost complex structure} $\JJC$ on $\TTC M$ is an almost complex structure on $\TTC M$ which is also orthogonal in the inner product $\<,\>$. Furthermore, $\JJC$ is \emph{integrable} and is called an \emph{extended complex structure} if the $+i$-eigensubbundle $L$ of $\JJC$ is involutive with respect to the bracket $[,]$.
\end{defn}
\begin{remark} \label{sub:hitchin}
\rm{In \cite{Severa}, the notion of \emph{exact Courant algebroid} is used for the Courant algebroid $\TTC M$ here, since it fits into the exact sequence:
$$0 \to T^*M \to \TTC M \xto a TM \to 0.$$
Let $s : TM \to \TTC M$ be a splitting with isotropic image (with respect to $\<,\>$), then it defines the $3$-form $H_s = 2\<s(X), [s(Y), s(Z)]\>$, which is closed. Then $s$ identifies the Courant algebroid $\TTC M$ with $\TT M$ with the $H_s$-twisted Courant bracket. When such a splitting is chosen, we will use the notations $\TT M$, $\JJ$, etc, and say that we have the corresponding \emph{generalized} structures.
The set of splittings is a torsor over $\Omega^2(M)$. The class $[H_s] \in H^3(M)$ does not depend on the choice of splitting and is called the \emph{\v Severa characteristic class} of $\TTC M$.
}
\end{remark}
\begin{defn}\label{subg:extcplx}
The subgroup of generalized symmetries which preserve an extended complex structure $\JJC$ is $\GGS_\JJC$.
Choose an (isotropic) splitting of $\TTC M$ and let $H$ be the corresponding twisting form, then we have $\GGS_\JJ \subset \GGS_H$.
\end{defn}

\subsection{Action on spinors} As shown in \cite{Gualtieri}, each maximally isotropic subbundle $L \subset \TT_\CC M$ corresponds to a pure line subbundle $U$ of the spin bundle $\wedge^\bullet T^*_\CC M$:
$$U = \Ann_C(L) := \{\rho \in \wedge^\bullet T^*_\CC M| \mf X \cdot \rho = \iota_X \rho + \xi \wedge \rho = 0 \text{ for all } \mf X = X + \xi \in L\},$$
where $\cdot$ stands for the Clifford multiplication.
A (nowhere vanishing) local section $\rho$ of $U$ is called a \emph{pure spinor} associated to the subbundle $L$. The integrability of $L$ with respect to the $H$-twisted Courant bracket is equivalent to the condition $d_H(\Gamma(U)) \subset \Gamma(U_1)$, where $U_1 = \Gamma(\TT_\CC M) \cdot U$ via Clifford multiplication. More explicitly, there is a unique local section $\mf Y = Y + \eta$ of $\bar L$, so that 
\begin{equation}\label{subg:integ}
d_H \rho = d\rho - H \wedge \rho = \mf Y \cdot \rho = \iota_Y \rho + \eta \wedge \rho,
\end{equation}
where we use the same convention for $d_H$ as that in \cite{Mathai1}. For a generalized complex structure $\JJ$, we say that $\rho$ \emph{defines} $\JJ$ iff $\rho$ is the pure spinor defining the $+i$-eigenbundle $L$ of $\JJ$.
\begin{lemma}
Let $\rho$ be a spinor, $\tilde \lambda = (\lambda, \alpha) \in \GGS$ and define $\tilde \lambda \circ \rho := (\lambda^{-1})^*(e^{-\alpha}\rho)$. Correspondingly, the infinitesimal action of $(X, A) \in \XXS$ on $\rho$ is 
\begin{equation}\label{gency:infinitesimal}
(X, A) \circ \rho = -\LLC_X \rho - A\wedge \rho.
\end{equation} Suppose that $\rho$ defines a generalized complex structure $\JJ$, then $\tilde \lambda \circ \rho$ defines the generalized complex structure $\tilde \lambda \circ \JJ$. We also have
$d_{\tilde \lambda \circ H} \tilde \lambda \circ \rho = \tilde \lambda \circ d_H \rho$.
\end{lemma}
{\it Proof:}
Straightforward computation from the definitions.
\qed

\section{Hamiltonian symmetries}\label{ham}

\subsection{Infinitesimal action}\label{ham:infinitesimal} As shown by Gualtieri \cite{Gualtieri}, the infinitesimal Hamiltonian fields on a generalized complex manifold $M$ can be defined by a complex valued function $F: M \to \CC$ as $X+\xi = \frac{1}{2}(d_L F + \bar{d_L F})$, where $d_L$ is the Lie algebroid differential defined on $L$ (see \eqref{app:liediff}), and the infinitesimal symmetry is generated by $(X, d\xi)$. 
It's straight forward to check that 
in the decomposition $\TT_\CC M = L \dsum L^*$, $d_L F = (dF)_{L^*} = dF + i\JJ (dF)$. Write $F = \Re F + i \Im F$, then
$$X+\xi = \frac{1}{2}(d_L F + \bar{d_L F}) =  \frac{1}{2}(dF + i\JJ (dF) + d\bar F - i\JJ (d\bar F)) = d(\Re F) - \JJ (d \Im F).$$
Let $X_F + \xi_F =  - \JJ (d \Im F)$, then the infinitesimal symmetry is $(X, d\xi) = (X_F, d\xi_F)$ when $H = 0$, in which there is no contribution from $\Re F$. For general $H$, we take into account of the embedding $\psi_H$ and obtain
\begin{defn}\label{ham:field}
Let $\mf X = X + \xi \in \Gamma(\TT M)$. The generalized symmetry preserving the $H$-twisted Courant bracket $[,]_H$ generated by $\mf X$ is $e^{t\mf X} := e^{t\psi_H(X, d\xi)}$.
Let $f: M \to \RR$ be a smooth function on an $H$-twisted generalized complex manifold $(M, \JJ)$. The \emph{Hamiltonian field} $\mf X_f$ generated by $f$ is $\JJ(df) = X_f + \xi_f$
and $f$ is a \emph{Hamiltonian function} defining the \emph{generalized Hamiltonian symmetry} $e^{t\mf X_f}$.
\end{defn}
\begin{remark}\label{ham:invariant}
\rm{Similar arguments as those in Chapter $5$ of \cite{Gualtieri} then imply that
$\{e^{t\mf X_f} \} \subset \GGS_\JJ$.
}
\end{remark}
\begin{lemma}\label{ham:prop} For $f, g: M \to \RR$, let $\{f, g\} = X_f(g)$. Then $\{,\}$ is a Poisson bracket and $(M, \JJ)$ is canonically a Poisson manifold. We have
\begineq\label{ham:propeq}
\iota_{X_f} df = \iota_{X_f} \xi_f = 0 \text{ and } [(X_f, d\xi_f), (X_g, d\xi_g)]_H = (X_{\{f,g\}}, d\xi_{\{f,g\}}).
\endeq
\end{lemma}
{\it Proof:} Since $idf +\JJ (df) = idf + X_f + \xi_f \in \Gamma(L)$, we find that $0 = \< idf + X_f + \xi_f,  idf + X_f + \xi_f\> = 2\iota_{X_f}(\xi_f+ idf)$. Separating the real and imaginary parts and we have the first equation in \eqref{ham:propeq}.

To show that $\{,\}$ is Poisson bracket, we compute the $H$-twisted Courant bracket of sections of the form $X_f + \xi_f + idf$ in $L$, without assuming $X_f \in \Gamma_H(TM)$:
\begin{equation*}
\begin{split}
& [X_f + \xi_f + idf, X_g + \xi_g + idg]_H \\ 
= & [X_f, X_g] + \LLC_{X_f}(\xi_g + i dg) - \iota_{X_g}d(\xi_f + i df) + \iota_{X_g}\iota_{X_f} H\\
= & [X_f, X_g] + \LLC_{X_f} \xi_g - \iota_{X_g} d\xi_f + \iota_{X_g}\iota_{X_f} H + i d \{f, g\}
\end{split}
\end{equation*}
It follows that $[X_f, X_g] = X_{\{f, g\}}$ and $\xi_{\{f, g\}} = \LLC_{X_f}\xi_g - \iota_{X_g} d\xi_f + \iota_{X_g}\iota_{X_f} H$, proving the second equation in \eqref{ham:propeq}. The rest of the lemma then follows easily from these.
\qed

\begin{remark}\label{ham:extend}
\rm{As mentioned in the introduction, several authors (e.g. \cite{Gualtieri}, \cite{Abouzaid} and \cite{Crainic}) have noticed that $(M, \{,\})$ is Poisson. The non-twisted version of the definition \ref{ham:field} and lemma \ref{ham:prop} also appeared in various forms in the works cited above. The Poisson structure does not depend on the splitting of $\TTC M$ since $X_f = a(\JJC (df))$.
}
\end{remark}

\subsection{Group of Hamiltonian symmetries}
Let $F_t$ (resp. $G_t$) be time-dependent Hamiltonian function on $(M, \JJ)$, with $\JJ(dF_t) = X_t + \xi_t$ (resp. $\JJ(dG_t) = Y_t + \eta_t$), which generates infinitesimal symmetries $\psi_H(X_t, d\xi_t)$ (resp. $\psi_H(Y_t, d\eta_t)$) and  path of generalized symmetries $(\lambda_t, \alpha_t)$ (resp. $(\mu_t, \beta_t)$). Then $(\lambda_t, \alpha_t)$ and $(\mu_t, \beta_t)$ are paths in $\GGS_\JJ$ (cf. remark \ref{ham:invariant}). 
\begin{lemma}\label{ham:compfunc}
The path of generalized symmetries $(\tau_t, \gamma_t)
= (\lambda_t\mu_t, \mu_t^*\alpha_t + \beta_t)$ is generated by $K_t = F_t + (\lambda_t^{-1})^*G_t$.
\end{lemma}
{\it Proof:} 
Let $(X_t, A_t) = (X_t, d\xi_t - \iota_{X_t} H)$ and $(Y_t, B_t) = (Y_t, d\eta_t - \iota_{Y_t} H)$. By \eqref{ham:time}, we have $\dot \alpha_t = \lambda_t^*A_t$ and $\dot\beta_t = \mu_t^*B_t$.
Consider $(Z_t, C_t)$ where $\dot\gamma_t = \rho_t^*C_t$ and $Z_t = X_t + \lambda_{t*} Y_t$ generates $\rho_t = \lambda_t\mu_t$. We only need to show that $C_t = d\zeta_t - \iota_{Z_t}H$ with $\JJ(dK_t) = Z_t + \zeta_t$. We first compute
$$\dot\gamma_t = \left.\frac{d}{ds}\right|_{s = 0}(\mu_t^*(\mu_t^{-1})^*\mu_{t+s}^*\alpha_{t+s} + \beta_{t+s})
=\mu_t^* \LLC_{Y_t} (\alpha_t) + \mu_t^*\dot\alpha_t + \dot \beta_t.$$
It follows that 
$$C_t = (\rho_t^{-1})^* \dot \gamma_t = (\lambda_t^{-1})^*\LLC_{Y_t} \alpha_t + A_t + (\lambda_t^{-1})^*B_t.$$
By definition, we have $X_t+\xi_t + idF_t \in \Gamma(L)$ as well as $Y_t + \eta_t + i dG_t \in \Gamma(L)$, then 
$$(\lambda_t, \alpha_t)\circ(Y_t + \eta_t + i dG_t) = \lambda_{t*}Y_t + (\lambda_t^{-1})^*(\eta_t + \iota_{Y_t}\alpha_t) + id((\lambda_t^{-1})^*G_t) \in \Gamma(L).$$
It follows that $\zeta_t = \xi_t + (\lambda_t^{-1})^*(\eta_t + \iota_{Y_t}\alpha_t)$. The lemma then follows from checking 
$$C_t = d(\xi_t + (\lambda_t^{-1})^*(\eta_t + \iota_{Y_t}\alpha_t)) - \iota_{X_t + \lambda_{t*}Y_t} H.$$
\qed

Completely parallel to symplectic case, we have the following:
\begin{theoremdefn}
\label{ham:gp}
Let $(M, \JJ)$ be $H$-twisted generalized complex manifold. 
The \emph{group of generalized Hamiltonian symmetries} $\Ham(M, \JJ) \subset \GGS_\JJ$ is the set of time-$1$ generalized symmetries generated by time-dependent Hamiltonian functions. 
\qed
\end{theoremdefn}

\subsection{Hamiltonian action}
We assume in the following that the Lie group $G$ is connected.
\begin{defn}\label{ham:moment}
The action of Lie group $G$ on $H$-twisted generalized complex manifold $(M, \JJ)$ is given by a group homomorphism to $\GGS_\JJ$. It is \emph{Hamiltonian} with \emph{moment map} $\mu : M \to \mf g^*$ 
if the induced geometric action on the Poisson manifold $(M, \{, \})$ is Hamiltonian with equivariant moment map $\mu$, 
so that the $G$-action is generated by $\psi_H(X_\mu, d\xi_\mu)$
, where $\JJ(d\mu) = X_\mu + \xi_\mu$. 
\end{defn}

\begin{prop}\label{ham:bfield}
Let the $G$-action on $(M, \JJ)$ be Hamiltonian with moment map $\mu$. Let $\JJ'$ be $B$-transformed generalized complex structure where $B \in \Omega^2(M)$. Then the same $G$-action is Hamiltonian with respect to $\JJ'$ iff $\LLC_{X_\mu} B = 0$, for which case the moment maps coincide.
\end{prop}
{\it Proof:}
We only need to check the condition $d\xi'_\tau - \iota_{X_\tau}H' = d\xi_\tau - \iota_{X_\tau} H$ for all $\tau \in \mf g$, where $\xi_\tau' = \xi_\tau + \iota_{X_\tau}B$ and $H' = H-dB$. It's equivalent to $\LLC_{X_\mu} B = 0$.
\qed

The following is obvious:
\begin{prop}\label{ham:prod}
For $i = 1, 2$, let $(M_i, \JJ_i)$ be $H_i$-twisted generalized complex manifold, with Hamiltonian $G_i$-action whose moment map is $\mu_i : M_i \to \mf g_i^*$. Then $(M = M_1 \times M_2, \JJ_1\dsum\JJ_2)$ is $H_1 \dsum H_2$-twisted generalized complex manifold, with Hamiltonian $G = G_1\times G_2$-action whose moment map is $\mu_1\times\mu_2 : M \to \mf g^*$. When $G_1 = G_2$, the diagonal action is Hamiltonian with moment map $\mu' = \mu_1 + \mu_2$.
\qed
\end{prop}

\section{Reduction by Hamiltonian action}\label{quot}

\subsection{General construction}\label{quot:constr} Let $G$ be a compact Lie group. Let $M$ be $H$-twisted generalized complex manifold with a Hamiltonian $G$ action with moment map $\mu: M \to \mf g^*$. Let $M_0 = \mu^{-1}(0)$ and $\JJ(d\mu_\tau) = X_\tau + \xi_\tau$ for $\tau \in \mf g$.
A few words on the notations below. The subscript $_\mu$ means that the associated object is valued in $\mf g^*$ (except those for $\TT_\mu$ and such) so that pairing with $\tau \in \mf g$ gives the cooresponding object associated to $\tau$. An expression such as $\theta \wedge \xi_\mu$, with $\theta$ being $\mf g$-valued connection form, invokes also the pairing between $\mf g$ and $\mf g^*$. Another equivalent way of unwinding $\theta\wedge \xi_\mu$ is to choose dual basis $\{\tau_i\}$ and $\{\tau^*_i\}$ of $\mf g$ and $\mf g^*$ and express $\theta = \sum_i \theta_i \tau_i$ and $\xi_\mu = \sum_i \xi_i \tau^*_i$, then $\theta\wedge \xi_\mu = \sum_i \theta_i \wedge \xi_i$.
\begin{assumption}\label{quot:assum} We list the assumptions that we'll use:
\begin{enumerate}
\item[$(0)$] $G$ action is preserves a splitting, i.e. $d\xi_\tau - \iota_{X_\tau}H = 0$ for all $\tau \in \mf g$,
\item[$(1)$] $0$ is a regular value of $\mu$,
\item[$(2)$] \emph{(the geometric part of)} $G$ action is free on $M_0$,
\end{enumerate}
\end{assumption}
\begin{remark}\label{quot:special}
\rm{ We'll drop the \emph{(the geometric part of)} in condition $(2)$ and simply say that $G$ action is free on $M_0$. We note the similarity of the condition $(0)$ to the exactness for Hamiltonian action in symplectic case, where $d\mu = -\iota_{X_\mu} \omega$. They both imply existence of an equivariantly closed extension, which is $H- u \xi_\mu$ here and $\omega+u \mu$ in symplectic geometry.
}
\end{remark}

\begin{lemma}\label{quot:linear}
Let $\VV = V\dsum V^*$ and $\JJ_V$ be a linear generalized complex structure. Let $\Ann(\cdot)$ denote the annihilating space in the pairing $\<,\>$.
Given subspace $K \subset V^*$, we assume
\begin{enumerate}
\item $K + \JJ_V(K) \subset \Ann(K, \JJ_V(K))$ and 
\item $\JJ_V(K) \inter V^* = \{0\}$.
\end{enumerate}
Then, $\<,\>$ descends to $\<,\>_K$ and $\JJ_V$ descends to $\JJ_K$ on the sub-quotient $\VV_K = \frac{\Ann(K \dsum \JJ_V(K))}{K \dsum \JJ_V(K)}$ as generalized complex structure. $\VV_K$ fits into the exact sequence, which splits non-canonically:
$$0 \to W_K^* \xto{a_K^*} \VV_K \xto{a_K} W_K \to 0, \text{ with } W_K = \frac{\Ann_V(K)}{N},$$
where $a_K$ is induced from projection $a: \VV \to V$ and $N = a\circ\JJ_V(K)$. Moreover, if 
\begin{enumerate}
\item[(3)] $\<N, c\circ\JJ_V(K)\> = 0$, where $c$ is the projection $\VV \to V^*$,
\end{enumerate}
the splitting map $W_K \to \VV_K$ can be defined by an element of $\wedge^2 V^*$.
\end{lemma}
{\it Proof:}
The subspace $\Ann(K, \JJ_V(K))$ is obviously closed under $\JJ_V$. By the first assumption, $\VV_K$ is well defined. The pairing  $\<,\>$ descends to a nondegenerate pairing $\<,\>_K$ since $(K, \JJ_V(K))$ is the null-space of the restriction of $\<,\>$. 
Also $\JJ_V$ descends to $\JJ_{K}$ on $\VV_{K}$ and is again generalized complex with respect to $\<,\>_{K}$. With assumption $(2)$, we have direct sum $K \dsum \JJ_V(K)$. The projection $a$ induces 
$$a_{K} : \VV_{K} \to W_{K} = \frac{\Ann_V(K)}{N}, \text{ with } \ker a_{K} = \frac{\Ann_{V^*}(N) \dsum \JJ_V(K)}{K \dsum \JJ_V(K)} \simeq W_K.$$
The fact that $a_{K}$ is surjective follows from the second assumption.
Note that $\Ann_{V^*}(N) \dsum \JJ_V(K)$ is isotropic with respect to $\<,\>$ and it follows that $\ker a_K$ is maximally isotropic with respect to $\<, \>_K$.

Let $\{u_l\}_{l = 1}^k$ be a basis of $K$ and $\JJ_V(u_l) = v_l + g_l$, then $\{v_l\}$ is a basis of $N$. Choose $\{v_j^* \in V^*\}$ so that $\<v_j^*, v_l\> = \delta_{jl}$, then the sequence is split by the map $W_{K} \to \VV_{K}$ induced by
$\Ann_V(K) \to \Ann(K, \JJ_V(K)) : w \mapsto w - \sum_{j = 1}^k\<w, g_j\> v_j^*$.
If furthermore, $\<N, c\circ \JJ_V(K)\> = 0$, then we define $B = \sum_{j=1}^k g_j \wedge v_j^* \in \wedge^2V^*$ and the map $w \mapsto w - \iota_wB$ splits the sequence.
\qed

\begin{theorem}\label{quot:algebroid}
With the assumptions \ref{quot:assum} $(1)$ and $(2)$, there is a natural extended complex structure on the quotient $Q = M_0/G$. Moreover, the extension sequence splits when pulled back to $M_0$, with the choice of a connection form on $M_0$.
\end{theorem}
{\it Proof:}
We first carry out the linear algebra for the bundles $\TT M |_{M_0}$ and $df \in T^*M |_{M_0}$. Consider the bundles over $M_0$:
$$\TT_\mu M_0 = \frac{\Ann(d\mu, \JJ (d\mu))}{\left(d\mu, \JJ(d\mu)\right)} = \frac{\Ann(d\mu, X_\mu + \xi_\mu)}{\left(d\mu, X_\mu + \xi_\mu\right)} \text{ and } W_\mu M_0 = \frac{\Ann_{TM}(d\mu)}{(X_\mu)}.$$
The assumption $(2)$ of lemma \ref{quot:linear} is given by the assumption \ref{quot:assum}-$(2)$. To see that $\TT_\mu M_0$ is well defined, we note that $\<d\mu_\tau, d\mu_\omega\> = \<\JJ(d\mu_\tau), \JJ(d\mu_\omega)\> = 0$ for all $\tau, \omega \in \mf g$. Then $\<d\mu_\tau, \JJ(d\mu_\omega)\> = \iota_{X_\omega} d\mu_\tau = \mu_{[\omega, \tau]} = 0$ on $M_0$ gives assumption $(1)$ of lemma \ref{quot:linear}. Thus lemma \ref{quot:linear} implies that $\JJ$ and $\<,\>$ descends to $\TT_\mu M_0$ and gives almost generalized complex structure $\JJ_\mu$ with respect to $\<, \>_\mu$ on the bundle. We also have the exact sequence:
$$0 \to W^*_\mu M_0 \xto{a_\mu^*} \TT_\mu M_0 \xto{a_\mu} W_\mu M_0 \to 0.$$
The bundle $\TT_\mu M_0$ is $G$-equivariant since both $\Span\{d\mu, \JJ(d\mu)\}$ and $\Ann(d\mu, \JJ(d\mu))$ are $G$-equivariant subbundles of $\TT M$.
The bundle $\TT_\mu M_0$ with the structure $\left(\JJ_\mu, \<,\>_\mu\right)$ descends to the bundle $\TTC_\mu Q$ with the structure $\left(\JJC_\mu, \<,\>_\mu\right)$ on $Q$.
Obviously, $\JJC_\mu$ is both complex and symplectic with respect to $\<,\>_\mu$. Since $W_\mu M_0$ naturally identifies with $\pi^*(TQ)$, the exact sequence above descends:
$$0 \to T^*Q \xto{a_\mu^*} \TTC_\mu Q \xto{a_\mu} TQ \to 0.$$
The image of $T^*Q$, i.e. $\ker a_\mu$ is maximally isotropic in the induced pairing. To split the pull-back sequence on $M_0$, 
we choose a connection form $\theta$ on $M_0$ and define:
\begineq\label{quot:splitpb}
\Ann_{TM}(d\mu) \to \Ann(d\mu, \JJ(d\mu)) : Y \mapsto Y - (\iota_Y\xi_\mu) \theta.
\endeq

Let $X+ \xi$ and $X'+\xi'$ be invariant sections of $\Ann(d\mu, X_\mu+ \xi_\mu)$, i.e. the following vanishing is true (see \eqref{subg:twistinf}):
$$\LLC_{X_\tau}\xi -\iota_X d\xi_\tau + \iota_X\iota_{X_\tau} H = 0, [X_\tau, X] = 0, \iota_X d\mu = 0 \text{ and } \iota_X \xi_\tau + \iota_{X_\tau}\xi = 0$$
and ditto for $X' + \xi'$. Direct computation shows that $[X+\xi, X'+\xi']_H$ satisfies the above vanishing equations, i.e. is again invariant and in $\Ann(d\mu, X_\mu + \xi_\mu)$. It's easy to see that $[X+\xi, d\mu_\tau]_H = 0$ and we compute:
$$[X+\xi, X_\tau + \xi_\tau]_H
= [X, X_\tau] + \LLC_X \xi_\tau - \LLC_{X_\tau} \xi - \frac{1}{2}d\(\iota_X \xi_\tau - \iota_{X_\tau}\xi\) + \iota_{X_\tau}\iota_X H
 = 0.
$$
We point out that the computations in this paragraph only use the vanishing equations. Now \eqref{app:couranteq3} implies
$$[X+\xi, k d\mu_\tau + l (X_\tau + \xi_\tau)]_H = X(k) d\mu_\tau + X(l)(X_\tau + \xi_\tau) \text{ for } k, l \in C^\infty(M_0)$$
and it follows that the $H$-twisted Courant bracket $[,]_H$ descends to $\TTC_\mu Q$ as a Courant bracket $[,]_\mu$. Because $\JJ$ is integrable, $\JJC_\mu$ is integrable under the induced bracket on $\TTC_\mu Q$. Thus $(\TTC_\mu Q, a_\mu, \<,\>_\mu, [,]_\mu, \JJC_\mu)$ is an extended complex structure on $Q$.
\qed

We show the effect of a $B$-transformation.
\begin{lemma}\label{quot:bfield}
With the assumption \ref{quot:assum} $(1)$ and $(2)$, let $(M, \JJ_1)$ be the $B_1$-transformed generalized complex structure for $B_1 \in \Omega^2(M)^G$
, then with all other choices fixed to be the same, the extended complex structure $\JJC_1$ constructed on the reduction from $\JJ_1$ is a $b$-transformation of $\JJC$ for some $b \in \Omega^2(Q)$.
\end{lemma}
{\it Proof:}
By proposition \ref{ham:bfield}, we see that $G$ action on $(M, \JJ_1)$ is Hamiltonian with the same moment map.
Let $\theta$ be the chosen connection form and choose a basis $\{\tau_j\}$ of $\mf g$ and $X_j = X_{\tau_j}$, we define the horizontal part of $B_1$ as 
$$B' = \prod_j(1-\theta_j \wedge \iota_{X_j}) B_1 = B_1 - \sum_j \theta_j \wedge \iota_{X_j} B_1 + \sum_{j < k} \theta_j\wedge \theta_k \cdot \iota_{X_k} \iota_{X_j} B_1,$$
where we interpret $(1-\theta_j \wedge \iota_{X_j})$ as operators on $\Omega^2(M)$. Consider the $B'$-transformed structure $\JJ'$ and the corresponding $(\TTC' Q, \JJC')$. Since $\iota_{X_j}B' = 0$ for all $j$, $\Ann(d\mu, \JJ'(d\mu)) = \Ann(d\mu, \JJ(d\mu))$ and the bundles $\TTC' Q$ and $\TTC Q$ are identical. Since $\LLC_{X_k}$ commutes with $\sum_j\theta_j \wedge \iota_{X_j}$, we see that there is $b \in \Omega^2(Q)$ so that $\pi^* b = B'$, the bracket on the extended tangent bundles $\TTC' Q$ and $\TTC Q$ are related by $b$-transformation and $\JJC'$ is the $b$-transform of $\JJC$. Now $(\TTC_1 Q, \JJC_1)$ and $(\TTC'Q, \JJC')$ are isomorphic since $(B'-B_1)$-transformation identifies corresponding structures on $\Ann(d\mu, \JJ'(d\mu))$ and $\Ann(d\mu, \JJ_1(d\mu))$.
\qed
\begin{remark}\label{quot:nonexact}
\rm{
We point out that the form $b$ might not be exact even if $B_1$ were exact. Suppose that $G = S^1$ and $B_1 = d\tilde\alpha$ for some $\tilde\alpha\in \Omega^1(M)^G$, then $\tilde b = d\tilde\alpha - \theta\wedge \iota_{X_\mu} d\tilde\alpha = d(\tilde\alpha - \iota_{X_\mu}\tilde\alpha \theta) +  \iota_{X_\mu}\tilde\alpha d\theta$, which descends to $b = d\alpha + f\Omega$ for some $\alpha \in \Omega^1(Q)$ and $f \in C^\infty(Q)$, where $\Omega$ is the curvature form of $\theta$.
}
\end{remark}

\subsection{Induced splitting}\label{quot:gencplx} Let $B = \theta\wedge \xi_\mu$ be the $2$-form defined locally near $M_0$, where $\theta$ is a connection form extending the one in the proof of theorem \ref{quot:algebroid}. 
Let
\begin{equation}\label{quot:horizonb}
B_1 = B - \frac{1}{2}\sum_{j<k} \theta_j \wedge \theta_k \cdot \iota_{X_k}\iota_{X_j} B = \theta\wedge\xi_\mu - \frac{1}{2}\sum_{j,k}\theta_j\wedge\theta_k \cdot \iota_{X_k}\xi_j.
\end{equation}
The last equality we used the fact that $\iota_{X_k} \xi_j + \iota_{X_j} \xi_k = \<\JJ(d\mu_j), \JJ(d\mu_k)\> = 0$. We compute
\begin{equation*}
\begin{split}
\iota_{X_l} B_1 & = \xi_l - \sum_{j}\theta_j \cdot \iota_{X_l}\xi_j - \frac{1}{2}\iota_{X_l}\sum_{j, k} \theta_j \wedge \theta_k \cdot \iota_{X_k}\xi_j \\
& = \xi_l - \sum_{j}\theta_j \cdot \iota_{X_l}\xi_j - \frac{1}{2}\(\sum_k \theta_k \cdot \iota_{X_k} \xi_l - \sum_j \theta_j \cdot \iota_{X_l} \xi_j\) = \xi_l.
\end{split}
\end{equation*} 
Apply the inverse transformation $e^{-B_1}$, we see that the generalized complex structure $\JJ_M = \JJ$ obtains extra twisting and becomes $H+dB_1$-twisted.
We have
\begin{corollary}\label{quot:descend}
With assumption \ref{quot:assum}, the extended complex structure $(\TTC_\mu Q, \JJC_\mu)$ splits and the locally defined form $H+dB_1$ as above descends to the quotient and gives $h$-twisted generalized complex structure on $Q$, where $H+dB_1 = \pi^*h$. The cohomology class $[h] \in H^3(Q)$ is independent of choice of connection. Let $B_2 \in \Omega^2(M)^G$ then $[h]$ is independent of $B_2$-transformation on $M$ as well.
\end{corollary}
{\it Proof:} Assumption $4.1$-$(0)$ implies that $\LLC_{X_\tau} B = 0$. Since $\LLC_{X_\tau}$ commutes with $\sum_k\theta_k\wedge \iota_{X_k}$, we see that $\LLC_{X_\tau} B_1 = 0$ as well. Let $\JJ_1 = e^{-B_1} \JJ e^{B_1}$, then $\Ann(d\mu, \JJ_1(d\mu)) = \Ann(d\mu, X_\mu)$ and the subquotient $\TT_{\mu, 1} M_0$ is canonically identified to $\pi^*\TT Q$, with structure $\JJC_{\mu, 1}$. The extension sequence naturally splits. Direct computation shows that the horizontal part of $B_1$ vanishes. By lemma \ref{quot:bfield}, $\TTC_\mu Q$ is isomorphic to $\TTC_{\mu, 1}Q$ and the structure $\JJC_\mu$ on $\TTC_\mu Q$ is identified with $\JJC_{\mu, 1}$ on $\TTC_{\mu, 1}Q$ by the bundle isomorphism induced from $e^{-B_1}$. Since $-B_1$-transform is orthogonal we see that $(\TTC_\mu Q, \JJC_\mu)$ splits.
To get the twisting of reduced structure, because $d(H + dB_1) = 0$, we only have to compute
$$\iota_{X_\tau}(H + dB_1)  = \iota_{X_\tau}H + \LLC_{X_\tau}B_1 - d\iota_{X_\tau}B_1 = \iota_{X_\tau} H - d\xi_\tau = 0.
$$
Thus $H+dB_1$ descends to $Q$, i.e. $H+dB_1 = \pi^*h$ for some $h \in \Omega_0^3(Q)$. The twisted Courant brackets obviously correspond.

For a different connection $\theta'$, let $B'_1$ be the corresponding form as in \eqref{quot:horizonb}, then the computation above shows that $\LLC_{X_\tau}(B'_1-B_1) = 0$ as well as $\iota_{X_\tau}(B'_1 - B_1) = 0$, thus there is $b \in \Omega^2(M)$ so that $\pi^* b = B_1'-B_1$.
Then the structures $\JJ_1$ and $\JJ'_1$ are related by $B_1'-B_1$-transformation. It follows that the reduced structures are related by $b$-transformation and the cohomology classes of the twistings are the same.
\qed
\begin{remark}\label{quot:coadjoint}
\rm{As with symplectic reduction, the reduction over coadjoint orbits can be performed via the shifting trick. Let $O_{\tau^*}$ be the coadjoint orbit passing through $\tau^* \in \mf g^*$, the condition \ref{quot:assum} is then stated for $\mf g^*_{\tau^*}$, where $G_{\tau^*}$ is the isotropy subgroup of the coadjoint action at $\tau^*$.
}
\end{remark}

\subsection{Generalized Calabi-Yau manifold} \label{quot:gency}
We show that if the manifold is generalized Calabi-Yau and the Hamiltonian action factors through $\Diff(M)$, then the reduced manifold is naturally generalized Calabi-Yau as well. 
\begin{corollary}\label{quot:gencyred}
Let $(M, \JJ)$ be an $H$-twisted generalized Calabi-Yau manifold. With the assumption \ref{quot:assum}, the reduced structure on $Q$ is an $h$-twisted Calabi-Yau manifold, where $h$ is as in corollary \ref{quot:descend}.
\end{corollary}
{\it Proof: } Choose a $d_H$-closed non-vanishing section of the pure spinor line bundle $U$. It follows from \eqref{gency:infinitesimal} that $\rho$ is invariant under the $G$-action. By Hamiltonian-ness, we see that
$$\iota_{X_\mu} \rho|_{M_0} = -i(d\mu \wedge \rho)|_{M_0} = 0$$
It follows that $\rho|_{M_0}$ descends to $\rho_\mu$ on $Q$. Let $B_1$ be as in corollary \ref{quot:descend} and $\tilde H = H + dB_1$, $\tilde \rho = e^{B_1}\rho$, then $d_{\tilde H} \tilde \rho = e^{B_1} d_H \rho = 0$. Since $B_1$-transformation doesn't change Hamiltonian-ness, it follows that $\tilde \rho$ descends to $\rho_\mu$ on $Q$ and $d_h \rho_\mu = 0$. That the $h$-twisted generalized complex structure on $Q$ is defined by $\rho_\mu$ and that $\rho_\mu$ is nowhere vanishing are then straight forward.
\qed

\subsection{Torus action}\label{quot:torus}
When $G = T$ is commutative, i.e. torus, we show that a Duistermaat-Heckman type formula is true when the action factors through $\Diff(M)$.
\begin{corollary}\label{quot:DHform}
With the assumption \ref{quot:assum}, for $G = T$ a torus, we have $h = h_0 + \Omega \wedge \zeta_\mu$, where $\tilde h_0 = \pi^*(h_0)$ is the horizontal part of $H$, $\tilde\zeta_\mu = \pi^*(\zeta_\mu)$ is the horizontal part of $\xi_\mu$ and $\Omega$ is the curvature form of the principle $T$-bundle over $Q$.
\end{corollary}
{\it Proof:} We compute:
\begin{equation*}\begin{split}
H+dB_1 = & H + d(\theta \wedge \xi_\mu - \frac{1}{2}\sum_{j,k}\theta_j\wedge\theta_k \cdot \iota_{X_k}\xi_j) \\
= & (H - \theta \wedge d\xi_\mu) + \pi^*\Omega \wedge \xi_\mu - \\
& - \frac{1}{2}\sum_{j,k}\(d\theta_j \wedge \theta_k \cdot \iota_{X_k}\xi_j - \theta_j \wedge d\theta_k \cdot \iota_{X_k}\xi_j + \theta_j\wedge \theta_k \cdot d\iota_{X_k}\xi_j\)\\
= & (H - \theta \wedge \iota_{X_\mu}H + \sum_{j<k} \theta_j \wedge \theta_k \cdot \iota_{X_k}\iota_{X_j} H) + \pi^*\Omega \wedge (\xi_\mu - \sum_k \theta_k \wedge \iota_{X_k}\xi_\mu) \\
=_\dag & \tilde h_0 + \pi^*\Omega \wedge \tilde \zeta_\mu.
\end{split}
\end{equation*}
The first parts on the two sides of equality $\dag$ are equal since $\LLC_{X_\tau} \xi_\mu = 0$ and abelian-ness imply that $\iota_{X_j}\iota_{X_k}\iota_{X_l} H = \iota_{X_j}\iota_{X_k} d\xi_l = 0$.
\qed

\subsection{Example: $\CC^2 \setminus\{(0, 0)\}$}\label{quot:exple}
We give an example of Hamiltonian $S^1$-action on a twisted generalized complex (in fact, it's generalized K\"ahler) manifold. As shown in Gualtieri \cite{Gualtieri}, the Hopf surface $S$ does not admit any generalized K\"ahler structure without twisting. In \cite{Gualtieri}, it's also shown that a twisted generalized K\"ahler structure can be given on $S$. We will put two $S^1$-actions on $M = \CC^2 \setminus \{(0,0)\}$ with the generalized K\"ahler structure lifted from $S$, so that the actions factor through $\Diff(M)$ and are Hamiltonian with respect to one of the twisted generalized complex structures.

Recall that Hopf surface $S = M / \ZZ$, where $\ZZ$ acts by scaling $n\circ z = 2^n z$. Let the metric on $M$ be $g = r^{-2}\tilde g$ where $\tilde g$ is the standard metric on $\CC^2$ and $r^2 = |z_1|^2 + |z_2|^2$, then it descends to $S$. We'll not say anything more about Hopf surface since we will mainly work on $M$. For more detail on how to get the following twisted generalized K\"ahler structure on $S$ (and thus on $M$), please consult \cite{Gualtieri}.

Let $J = \(\begin{smallmatrix}0 & -1 \\ 1 & 0 \end{smallmatrix}\)$ and $z = (z_1, z_2) = (x_1 + i y_1, x_2 + iy_2) = (x_1, y_1, x_2, y_2)$ be the coordinates in $\CC^2$, then we may write down the two generalized complex structures on $M$ as following:
\begin{equation*}
\JJ_1 =  \(\begin{matrix}
	0 & 0 & r^2 J & 0 \\
	0 & -J & 0 & 0 \\
	r^{-2} J & 0 & 0 & 0 \\
	0 & 0 & 0 & -J
	\end{matrix}\),
\JJ_2 = \(\begin{matrix}
	J & 0 & 0 & 0 \\
	0 & 0 & 0 & -r^2 J \\
	0 & 0 & J & 0 \\
	0 & -r^{-2} J & 0 & 0
	\end{matrix}\),
\end{equation*}
where the labelling on rows are $(T z_1, T z_2, T^*z_1, T^*z_2)^T$. They are both $H$-twisted generalized complex structures, where
\begin{equation*}
H = \frac{2}{r^4}(y_1dx_1 - x_1dy_1 + y_2dx_2 - x_2dy_2)\wedge(dx_1\wedge dy_1 + dx_2\wedge dy_2).
\end{equation*}
It's shown in \cite{Gualtieri} that $H$ represents nontrivial cohomology in $M$. In fact, let $z_1 = r e^{i\phi_1}\sin \lambda$ and $z_2 = r e^{i\phi_2}\cos \lambda$, where $\lambda \in [0, \frac{\pi}{2}]$ and $\phi_j \in [0, 2 \pi)$, then we compute that $H = -\sin (2\lambda) d\lambda \wedge d\phi_1 \wedge d\phi_2$ on $M$. 

Now let $f = \ln r$, then
$$df = \frac{1}{r^2}(x_1dx_1 + y_1 dy_1 + x_2dx_2 + y_2dy_2),$$
and it's not hard to check that $\JJ_1(df) = X_1 +\xi_1$ where
$$X_1 = x_1 \frac{\partial}{\partial y_1} - y_1\frac{\partial}{\partial x_1} \text{ and } \xi_1 = \frac{1}{r^2}(y_2dx_2 - x_2dy_2).$$
Then $X_1$ generates the action of $S^1$ rotating the $z_1$ plane. Direct computation then shows that $\iota_{X_1}H = d\xi_1$, i.e. the action factors through $\Diff(M)$. In fact, the same map $f$ with respect to $\JJ_2$ gives another Hamiltonian $S^1$-action, which rotates $z_2$ plane in the negative direction. We write down the components of $\JJ_2(df) = X_2 + \xi_2$:
$$X_2 = -x_2 \frac{\partial}{\partial y_2} + y_2\frac{\partial}{\partial x_2} \text{ and } \xi_2 = - \frac{1}{r^2}(y_1dx_1 - x_1dy_1).$$

Note that the action of $X_j$ fixes $\{z_j = 0\}$, where $j = 1, 2$. Thus the actions on the level set $f^{-1}(\ln r)$ always has a fixed $S^1$. We may consider reduction on $M_1 = \CC^2 \setminus \{z_1 = 0\}$, for example, with respect to the action generated by $X_1$. First of all, topologically, the quotient $Q_r$ is an open disc in $\CC$ of radius $r$. Thus the generalized complex structure on it can't be twisted due to dimension reasons. In fact, in coordinates $(\lambda, \phi_1, \phi_2)$, we have on $\{|z| = r, z_1 \neq 0\}$:
$$\xi_1 = - \cos^2 \lambda d\phi_2 \text{ and } H = -dB, \text{ where } B = -\cos^2 \lambda d\phi_1 \wedge d\phi_2 = d\phi_1 \wedge \xi_1.$$
Noting that the connection form on $\{|z| = r, z_1 \neq 0\}$ can be chosen as $d\phi_1$, we find that the locally defined form $H+dB$ is indeed $0$. As our computation of $H$ shows that $H = -dB$ is true in the whole $M_1$, the generalized complex structure we are reducing is in fact $B$-transformed from an untwisted one, say $\JJ_1' = e^{-B} \JJ_1 e^B$. In matrices, we have
\begin{equation*}\begin{split}
B = \(\begin{matrix} 0 & b \\ -b^T & 0 \end{matrix}\), & \text{ where } b = \frac{1}{r^2 |z_1|^2} \(\begin{matrix} -y_1 \\ x_1 \end{matrix}\)\(\begin{matrix} -y_2 & x_2\end{matrix}\) \\
& \text{ and rows of } B \text{ is labelled as } \begin{matrix} T^*z_1 \\ T^*z_2 \end{matrix}, \text{ and }
\end{split}\end{equation*}
$$\JJ_1' = \(\begin{matrix} 
	0 & r^2Jb & r^2 J & 0 \\ 
	0 & \boxed{-J} & 0 & \boxed{0} \\
	r^{-2}J & b J & 0 & 0 \\
	J b^T & \boxed{0} & r^2 b^T J & \boxed{-J}
	\end{matrix}\),
\JJ_2' = \(\begin{matrix}
	J & 0 & 0 & 0 \\
	r^2 J b^T & \boxed{0} & 0 & \boxed{-r^2 J} \\
	0 & 0 & J & r^2 b J \\
	b^T J & \boxed{-r^{-2} J} & 0 & \boxed{0} \\
	\end{matrix}\)
$$
It's then easy to check that $\JJ_1'(df) = X_1$ and that the quotient structure on $Q_r$, which is given by the boxed terms in $\JJ_1'$ above, is the restriction of the opposite complex structure on $\CC$. In fact, we also see that $\JJ_2'$ descends to quotient as well, becomes $r^{-2}$ times the restriction of the opposite symplectic structure on $\CC$. In other words, the quotient $Q_r$ can be identified as the open unit disc $D \subset \CC$ with the opposite K\"ahler structure. 

Similarly one may show that the reduction with respect to the action generated by $X_2$ gives open unit disc with the induced K\"ahler structure.

\section{Generalized complex cutting}\label{cut}
We describe here the cutting construction using the reduction construction in the previous section. We only state it for $S^1$-actions and $\CC$, while cutting with torus group and toric varieties can be constructed similarly as in the classical symplectic case (cf. \cite{Meinrenken}, \cite{Lerman}). Moreover, to simplify statement, we are going to insist on the full assumption \ref{quot:assum} (cf \S\ref{quot:gencplx}). 

Let $(M, \JJ)$ be an $H$-twisted generalized complex manifold with Hamiltonian $S^1$-action, whose moment map is $f$. Let $(\CC, \omega)$ be the symplectic manifold with standard symplectic structure and $S^1$-action whose moment map is denoted $g$. By shifting the images, we may assume that the moment maps $f$ and $g$ have ranges $[0, a]$ and $[0,\infty)$ for some $a > 0$. Now consider $M \times \CC$ with the product generalized complex structure $\tilde\JJ$ and the moment map $F = f+g : M \times \CC \to [0, \infty)$. Then $\tilde \JJ$ is $\pi_1^*H$-twisted, where $\pi_1$ is the projection to the first factor. Suppose that $\eps$ is a regular value of $f$ and $S^1$ acts freely on $f^{-1}(\eps)$, then theorem \ref{quot:descend} provides a twisted generalized complex structure $\JJ_\eps$ on $\bar {M^-_\eps} = f^{-1}([0, \eps]) / \sim$ where $\sim$ is the equivalent relation on the boundary given by $S^1$-action. As sanity check, we first show:
\begin{lemma}\label{cut:trivial}
When $a < \eps$ we have $f^{-1}(\eps) = \emptyset$ and $\bar{M^-_\eps} = M$ naturally as twisted generalized complex manifold.
\end{lemma}
{\it Proof:}
Let $\JJ(df) = X_f + \xi_f$ and $\JJ_\omega(dg) = Y_g$, then $\tilde \JJ(dF) = X_f\dsum Y_g + \xi_f \dsum 0$. Consider the following map:
$$\Phi: M\times I \times S^1 \to M \times \CC^\times : (m, a, \lambda) \mapsto (\phi(\lambda, m), a - f(m), \lambda),$$
which identifies $M\times I \times S^1$ to a neighbourhood of the level set $F^{-1}(\eps)$, where $\phi(\lambda, m)$ is the action of $\lambda\in S^1$ on $m$ and the coordinate on $\CC$ is given by $z = e^{2\pi i\lambda} \sqrt{2a}$. Let the $S^1$-action on the domain be trivial on the first and second factors and the multiplication on the third. Then $\Phi$ is equivariant. Furthermore, $\Phi\circ F = \pi_2$ the projection to the second factor. Pull everything back to $M\times I \times S^1$, then topologically the quotient is simply the projection to the first factor. Direct computation gives $(\Phi^*\tilde\JJ) (d\pi_2) = \frac{\partial}{\partial \lambda}$. The horizontal part of $\Phi^*\tilde\JJ$, i.e. restriction to any $M \times \{(a, \lambda)\}$, is identical to $\JJ$ due to invariance of $\JJ$ under $S^1$ action. Thus $\bar{M_\eps^-} = M$ as twisted generalized complex manifold.

To get the twisting form, we can also let $d\lambda$ denote the trivial connection on $F^{-1}(\eps)$ and consider the form $H' = H\dsum 0 + d(d\lambda \wedge (\xi_f \dsum 0))$. Then direct computation shows that $\Phi^*H' = \pi_1^* H$.
\qed

For the general case, we note that the principle $S^1$-bundle over $M_\eps^- = f^{-1}([0,\eps))$ is trivializable while the Chern class of the associated $\CC^1$-bundle of the principle $S^1$-bundle over $\bar{M_\eps^-}$ is the Poincar\'e dual of the reduced submanifold $Q_\eps = f^{-1}(\eps) / \sim$. 
Let $\eps' < \eps'' \in (0, \eps)$, so that $\eps - \eps'$ is small enough and $U_1=\bar{M_\eps^-} \setminus f^{-1}([0, \eps'])$ can be identified (equivariantly) with a neighbourhood of the normal bundle of $Q_\eps$ in $\bar{M_\eps^-}$, where the action is linear on the fiber. Let $U_2 = M_{\eps''}^-$, then $\{U_1, U_2\}$ is an open covering of $\bar{M_\eps^-}$. Let $\{\lambda_1, \lambda_2\}$ be a partition of unity subordinates to the covering $\{U_i\}$, so that $\lambda_i = 1$ in $U_i \setminus (U_1 \inter U_2)$. Let $\theta_\eps$ be the connection form on $f^{-1}(\eps)$ and define the reduction $Q_\eps$. Let $\theta_1$ be the pull back of $\theta_\eps$ via the bundle projection to $U_1$. Let $\theta_2$ be the trivial connection and $\theta$ be the connection form on $\bar{M_\eps^-}$ constructed from $\theta_i$ and the partition of unity $\lambda_i$.
\begin{prop}\label{cut:subfold}
With the above choice of connection form $\theta$, $Q_\eps$ is naturally a (twisted) generalized complex submanifold of $\bar{M_\eps^-}$ in the sense of \cite{Ben-Bassat}. In an open set away from $Q_\eps$, $\bar{M_\eps^-}$ retains the original generalized complex structure.
\end{prop}
{\it Proof:}
Since both $Q_\eps$ and $\bar{M_\eps^-}$ are already equipped with twisted generalized complex structures, we only have to show that they are indeed compatible. We show this by exhibiting $\TT Q$ as a natural subspace of $\TT \bar{M_\eps^-}$ at each point with the induced brackets and inner products. By construction, over $f^{-1}(\eps)$ there is natural inclusion $\TT_f M \subset \TT_F(M\times \CC)$ given by $Z+ \eta \mapsto Z \dsum 0 + \eta \dsum 0$. It induces inclusions $\Ann_Q(X_f + \xi_f, df) \subset \Ann_M(X_f \dsum Y_g + \xi_f \dsum 0, df \dsum dg)|_{f^{-1}(0)}$ and $(X_f + \xi_f, df) \subset (X_f \dsum Y_g + \xi_f \dsum 0, df \dsum dg)$, since $Y_g = 0$ and $dg = 0$ over $f^{-1}(\eps)$, where $\Ann_{\bullet}$ denotes the annihilator considered when constructing $Q$ or $\bar{M_\eps^-}$. All the structures restrict. By the choice of connection form, on $M_\eps$ we have natural inclusion of $\TT Q \subset \TT\bar{M_\eps^-}$ as generalized complex subspace in the sense of \cite{Ben-Bassat}.
\qed

\begin{remark}\label{cut:change}
\rm{ We note that the cutting construction does not need $S^1$-action on the whole manifold, instead, only a local $S^1$-action near the cut suffices, just as the symplectic case. That change of twisting does occur with the change of level sets is easy to see, when we note the two extreme cases, $M_\eps = \emptyset$ and $M_\eps = M$.
}
\end{remark}

\section{Appendix}\label{app}
\subsection{Courant and Lie algebroids}\label{app:broid} 
The material in this subsection is taken from 
\cite{Gualtieri}, where 
references can be found. Everything here can be complexified and get the corresponding complex notion.
\begin{defn}\label{app:liedefn}
A \emph{Lie algebroid} $L$ \emph{over manifold} $M$ is a vector bundle $\pi_L: L \to M$ with Lie bracket $[,]_L$ on the space of sections $\Gamma(L)$ as well as an \emph{anchor} homomorphism $a_L : L \to TM$ which satisfies the following:
\begin{eqnarray}
\label{app:lieeq1} a_L([X,Y]_L) & = & [a_L(X), a_L(Y)] \text{ for } X, Y \in \Gamma(L) \\
\label{app:lieeq2} [X, fY]_L & = & f\cdot [X, Y]_L + a_L(X)(f)\cdot Y \text{ for } f \in C^{\infty}(M).
\end{eqnarray}
\end{defn}
\begin{defn}\label{alg:liecomplex}
Let $L$ be a Lie algebroid over $M$ of rank $l$. Let $\Omega^k(L) = \Gamma(\wedge^kL^*)$, then the natural differential complex $(\Omega^{\bullet}(L), d_L)$ is defined as:
\begin{equation}\label{app:lieseq}
\Omega^{\bullet}(L) = \{0 \to \Omega^0(L) \xto{d_L} \Omega^1(L) \xto{d_L} \cdots \xto{d_L} \Omega^l(L) \to 0 \},
\end{equation}
where $d_L$ is given by the following algebraic formula for $\sigma \in \Omega^k(L)$ and $X_{0 \ldots k} \in \Gamma(L)$,
\begin{equation}\label{app:liediff}
\begin{split}
\(d_L\sigma\)(X_0, \ldots, X_k) = & \sum_{i}(-1)^i a_L(X_i) \sigma(X_0, \ldots, \hat{X}_i, \ldots, X_k) + \\
& + \sum_{i<j} (-1)^{i+j}\sigma([X_i, X_j]_L, X_0, \ldots, \hat{X}_i, \ldots, \hat{X}_j, \ldots, X_k).
\end{split}
\end{equation}
\end{defn}
\begin{defn}\label{app:courant}
A \emph{Courant algebroid} $E$ \emph{over manifold} $M$ is a vector bundle $\pi_E : E \to M$ with nondegenerate symmetric bilinear form $\<,\>_E$, a skew-symmetric bracket $[,]_E$ on $\Gamma(E)$ as well as an \emph{anchor} homomorphism $a_E : E \to TM$, which induces differential operator $\DDC : C^\infty(M) \to \Gamma(E)$ by $\<\DDC f, A\>_E = \frac{1}{2}a_E(A) f$ for all $f \in C^\infty(M)$ and $A \in \Gamma(E)$. They satisfy the following compatibility conditions for all $A, B, C \in \Gamma(E)$ and $f, g \in C^\infty(M)$:
\begin{eqnarray}
\label{app:couranteq1} a_E([A, B]_E) & = & [a_E(A), a_E(B)] \\
\label{app:couranteq2} \Jac(A, B, C) & = & \DDC (\Nij(A, B, C)) \\
\label{app:couranteq3} [A, fB]_E & = & f[A, B]_E + (a_E(A)f)B - \<A, B\>_E\DDC f \\
\label{app:couranteq4} a_E \circ \DDC & = & 0, \text{ i.e. } \<\DDC f, \DDC g\>_E = 0 \\
\label{app:couranteq5} a_E(A)\<B, C\>_E & = & \<[A, B]_E + \DDC\<A, B\>_E, C\>_E + \<B, [A, C]_E + \DDC \<A, C\>_E\>_E.
\end{eqnarray}
where
\begin{eqnarray}
\label{app:couranteq6} \Jac(A, B, C) & = & [[A, B]_E, C]_E + c.p. \\
\label{app:couranteq7} \Nij(A, B, C) & = & \frac{1}{3}\<[A, B]_E, C\>_E + c.p.
\end{eqnarray}
\end{defn}
In \cite{Uchino}, it's shown that with \eqref{app:couranteq1} and \eqref{app:couranteq5}, the Leibniz rule for $\DDC$, i.e. $\DDC(fg) = f\DDC(g) + \DDC(f) g$, implies the definition of $\DDC$ as well as \eqref{app:couranteq3} and \eqref{app:couranteq4}. An example of Courant algebroid is $\TT M$ with the  ($H$-twisted) Courant brackets and the natural pairing.

\subsection{Lie algebra extension}\label{app:coho}
The material in this subsection is taken from \cite{Knapp}, especially \S IV.$2$ and \S IV.$3$, where more detail is available. Although the statements in \cite{Knapp} are all for finite dimensional complex Lie algebras and their representations, they obviously adapt to our needs in the main text.

Let $\mf g$ be a finite dimensional Lie algebra and $V$ a representation of $\mf g$. The vector space $C^n(\mf g, V)$ of \emph{$n$-cochains} is
\begineq\label{app:cochain}
C^n(\mf g, V) = \Hom(\wedge^n \mf g, V),
\endeq
and the \emph{coboundary operator} $\partial_n : C^n(\mf g, V) \to C^{n+1}(\mf g, V)$ is defined by
\begineq\label{app:cobdry}
\begin{split}
(\partial_n \omega)&(Y_0 \wedge \ldots \wedge Y_n) = \sum_{i = 0}^n(-1)^l Y_i(\omega(Y_1 \wedge \ldots \wedge \hat Y_i \wedge \ldots \wedge Y_n)) + \\
& + \sum_{i<j}(-1)^{i+j} \omega([Y_i, Y_j]\wedge Y_1 \wedge \ldots \wedge \hat Y_i \wedge \ldots \wedge \hat Y_j \wedge \ldots \wedge Y_n).
\end{split}
\endeq
Then it's easy to check that $\partial_n\partial_{n-1} = 0$ and we then define
\begin{defn}\label{app:liecoho}
The space of \emph{$V$-valued cocycles (resp. coboundaries)} are $Z^n(\mf g, V) = \ker (\partial_n)$ (resp. $B^n(\mf g, V) = \img (\partial_{n-1})$). The \emph{$n$-th cohomology of $\mf g$ with coeffecients in $V$} is $H^n(\mf g, V) = Z^n(\mf g, V) / B^n(\mf g, V)$.
\end{defn}

Let $\mf g$ and $\mf a$ be Lie algebras, where $\mf a$ is assumed to be abelian, i.e. $[,] = 0$ in $\mf a$.
\begin{defn}\label{app:lieext}
A Lie algebra $\mf h$ is an \emph{extension} of $\mf g$ by $\mf a$ if there is an exact sequence of Lie algebras:
$$0 \to \mf a \xto{i} \mf h \xto{\rho} \mf g \to 0,$$
i.e. the maps $i$ and $\rho$ are Lie algebra homomorphisms and the image of $\mf a$ is an ideal in $\mf h$ which coincides with the kernel of $\rho$. Two extensions $\mf h$ and $\mf h'$ are \emph{equivalent} if there is a Lie algebra isomorphism $\sigma : \mf h \to \mf h'$, which together with the identity maps on the other terms gives isomorphism of the extension sequences.
\end{defn}

\begin{theorem}\label{app:lieextclass}
Let $\mf h$ be a Lie algebra extension of $\mf g$ by $\mf a$ with $\mf a$ abelian, then $\mf a$ is naturally a representation of $\mf g$ by $i(X\circ Y) = [\rho^{-1}(X), i(Y)]$, where $X \in \mf g$ and $Y \in \mf a$. The equivalent classes of extensions $\mf h$ which give rise to the same action of $\mf g$ on $\mf a$ are classified by $H^2(\mf g, \mf a)$.
\end{theorem}
More explicitly, for any cocycle $\omega \in Z^2(\mf g, \mf a)$ we define $\mf h = \mf g \dsum \mf a$ as vector space and the Lie bracket on $\mf h$ is given by:
\begineq\label{app:lieextbrac}
[(X, A), (Y, B)] = ([X, Y], X\circ B - Y \circ A + \omega(X, Y)).
\endeq


\begin{thebibliography}{90}
%
\bibitem{Abouzaid}
Mohammed~Abouzaid and Mitya~Boyarchenko,
Local structure of generalized complex manifolds,
{\it math.DG/0412084}, 2004
%
\bibitem{Albert}
Claude~Albert,
Le th\'eor\`eme de r\'eduction de Marsden-Weinstein en g\'eom\'etrie cosymplectique et de contact,
{\it J. Geom. Physics}, {\bf 6}(4):627--649, 1989.
%
\bibitem{Alekseev}
Anton~Alekseev, Anton~Malkin, and Eckhard~Meinrenken,
Lie group valued moment maps
{\it J. Differential Geom.} {\bf 48}(3):445--495, 1998.
%
\bibitem{Alekseev1}
Anton~Alekseev, Eckhard~Meinrenken and Christopher~Woodward,
Group-valued equivariant localization,
{\it Invent. Math.}, {\bf 140}(2):327--350, 2000.
%
\bibitem{Alekseev2}
Anton~Alekseev, Eckhard~Meinrenken and Christopher~Woodward,
 The Verlinde formulas as fixed point formulas, 
{\it J. Symplectic Geom.}, {\bf 1}(1):1--46, 2001.
%
\bibitem{Apostolov}
Vestislav~Apostolov and Simon~Salamon,
K\"ahler reduction of metrics with holonomy $G\sb 2$,
{\it Comm. Math. Phys.}, {\bf 246}(1):43-61, 2004.
%
\bibitem{Ben-Bassat}
Oren~Ben-Bassat and Mitya~Boyarchenko,
Submanifolds of generalized (almost) complex manifolds,
{\it J. Symplectic Geom.}, {\bf 2}(3):309--355, 2004.
%
\bibitem{Bielawski}
Roger~Bielawski,
Twistor quotients of hyperk\"ahler manifolds,
{\it Quaternionic structures in mathematics and physics (Rome, 1999)}, 7--21 (electronic),
Univ. Studi Roma ``La Sapienza'', Rome, 1999. 
%
\bibitem{Blankenstein1}
Guido~Blankenstein,
Implicit Hamiltonian systems: symmetry and interconnection,
{\it Ph.D. Thesis}, University of Twente, The Netherlands, November 2000.
%
\bibitem{Blankenstein2}
Guido~Blankenstein and Arjan~van~der~Schaft,
Symmetry and reduction in implicit generalized Hamiltonian systems,
{\it Rep. Math. Phys.}, {\bf 47}(1):57--100, 2001.
%
\bibitem{Blankenstein3}
Guido~Blankenstein and Tudor~S.~Ratiu,
Singular reduction of implicit Hamiltonian systems,
{\it Rep. Math. Phys}, {\bf 53}(2):211--260, 2004.
%
\bibitem{Bott}
Raoul~Bott, Susan~Tolman and Jonathan~Weitsman,
Surjectivity for Hamiltonian loop group spaces,
{\it Invent. Math.}, {\bf 155}(2):225--251, 2004.
%
\bibitem{Burns}
Daniel~Burns, Victor~Guillemin and Eugene~Lerman,
Kaehler cuts,
{\it math.DG/0212062}, 2002.
%
\bibitem{Bursztyn}
Henrique~Bursztyn, Gil~R.~Cavalcanti and Marco~Gualtieri,
Reduction of Courant algebroids and generalized complex structures, 
{\it math.DG/0509640}, 2005.
%
\bibitem{Crainic}
Marius~Crainic,
Generalized complex structures and Lie brackets,
{\it math.DG/0412097}, 2004.
%
\bibitem{Dragulete}
Oana~Dragulete and Liviu~Ornea,
Non-zero contact and Sasakian reduction,
{\it Differential Geom. Appl.}, {\bf 24}(3):260--270, 2006.
%
\bibitem{Geiges}
Hansj\"org~Geiges,
Constructions of contact manifolds,
{\it Math. Proc. Cambridge Philos. Soc.}, {\bf 121}(3):455--464, 1997.
%
\bibitem{Gini}
Rosa~Gini, Liviu~Ornea and Maurizio~Parton,
Locally conformal K\"ahler reduction.
{\it J. Reine Angew. Math.}, {\bf 581}(1):1--21, 2005.
%
\bibitem{Gini1}
Rosa~Gini, Liviu~Ornea, Maurizio~Parton and Paolo~Piccini,
Reduction of Vaisman structures in complex and Quaternionic geometry,
{\it math.DG/0502572}, 2005.
%
\bibitem{Ginzburg}
Viktor~L.~Ginzburg,
Some remarks on symplectic actions of compact groups,
{\it Math. Z.}, {\bf 210} (4):625--640, 1992.
%
\bibitem{Grantcharov}
Gueo~Grantcharov and Liviu~Ornea,
Reduction of Sasakian manifolds,
{\it J. Math. Physics}, {\bf 42}(8):3808--3816, 2001.
%
\bibitem{Gualtieri}
Marco~Gualtieri,
Generalized complex geometry,
Oxford thesis,
{\it math.DG/0401221}, 2004.
%
\bibitem{Hitchin}
Nigel~Hitchin,
Generalized Calabi-Yau manifolds, 
{\it Q. J. Math.}  {\bf 54} (3):281--308, 2003.
%
\bibitem{Hitchin2}
Nigel~Hitchin,
Instantons, Poisson structures and generalized K\"ahler geometry,
{\it Comm. Math. Phys.}, {\bf 265}(1):131--164, 2006.
%
\bibitem{Hitchin4}
Nigel~Hitchin, Anders~Karlhede, Ulf~Lindstr\"om and Martin~Ro{\v c}ek,
Hyper-K\"ahler metrics and supersymmetry.
{\it Comm. Math. Phys.} {\bf 108} (4):535--589, 1987.
%
\bibitem{Huebschmann2}
Johannes~Huebschmann,
Severi varieties and holomorphic nilpotent orbits,
{\it math.DG/0206143}, 2004.
%
\bibitem{Huebschmann}
Johannes~Huebschmann,
Stratified Kaehler structures on adjoint quotients,
{\it math.DG/0404141}, 2004.
%
\bibitem{Huebschmann0}
Johannes~Huebschmann,
K\"ahler spaces, nilpotent orbits, and singular reduction,
{\it Mem. Amer. Math. Soc.}, {\bf 172} (2004), no.~814, vi+96 pp.
%
\bibitem{Huebschmann1}
Johannes~Huebschmann, 
Singular Poisson-K\"ahler geometry of Scorza varieties and their secant varieties, 
{\it Differential Geom. Appl.}, {\bf 23} (1):79--83, 2005.
%
\bibitem{Kapustin}
Anton~Kapustin and Yi~Li,
Topological sigma-models with $H$-flux and twisted generalized complex manifolds,
{\it hep-th/0407249}, 2004.
%
\bibitem{Kirwan}
Frances~Kirwan,
Momentum maps and reduction in algebraic geometry,
{\it Differential Geom. Appl.} {\bf 9} (1-2):135--171, 1998.
%
\bibitem{Knapp}
Anthony~W.~Knapp, 
Lie groups, Lie algebras, and cohomology,
Princeton Univ. Press, Princeton, NJ, 1988.
%
\bibitem{Lerman1}
Eugene~Lerman,
Symplectic cuts,
{\it Math. Res. Lett.} {\bf 2} (3):247--258, 1995.
%
\bibitem{Lerman2}
Eugene~Lerman,
Contact cuts,
{\it Israel J. Math.} {\bf 124} 77--92, 2001.
%
\bibitem{Lerman}
Eugene~Lerman, Eckhard~Meinrenken, Susan~Tolman and Chris~Woodward,
Non-abelian convexity by symplectic cuts,
{\it Topology}, {\bf 37} (2):245--259, 1998. 
%
\bibitem{Lin}
Yi~Lin and Susan~Tolman,
Symmetries in generalized K\"ahler geometry,
{\it math.DG/0509069}, 2005.
\bibitem{Lindstrom1}
Ulf~Lindstr\"om and Martin~Ro{\v c}ek,
Scalar tensor duality and $N=1,\,2$ nonlinear $\sigma $-models.
Nuclear Phys. B 222 (1983), no. 2, 285--308.
%
\bibitem{Lindstrom}
Ulf~Lindstr\"om, Martin~Ro{\v c}ek and Rikard~von Unge,
Hyper-K\"ahler quotients and algebraic curves,
{\it J. High Energy Phys.}, {\bf 2000}(1) paper 22, 22 pages.
\bibitem{Lu}
Jiang-Hua~Lu,
Momentum mappings and reduction of Poisson actions,
In {\it Symplectic geometry, groupoids, and integrable systems (Berkeley, CA, 1989)},  209--226, 
Math. Sci. Res. Inst. Publ., 20, Springer, New York, 1991.
%
\bibitem{Marsden2}
Jerrold~E.~Marsden and Tudor~Ratiu,
Reduction of Poisson manifolds,
{\it Lett. Math. Phys.} {\bf 11}(2):161--169, 1986.
%
\bibitem{Marsden1}
Jerrold~E.~Marsden and Alan~Weinstein,
Reduction of symplectic manifolds with symmetry,
{\it Rep. Mathematical Phys.} {\bf 5}(1):121--130, 1974.
%
\bibitem{MarsdenY1}
Jerrold~E.~Marsden and Hiroaki~Yoshimura,
Dirac structures in Mechanics, Part I: Implicit Lagrangian systems,
{\it preprint}, 2005.
%
\bibitem{MarsdenY2}
Jerrold~E.~Marsden and Hiroaki~Yoshimura,
Dirac structures in Mechanics, Part II: Variational structures,
{\it preprint}, 2005.
%
\bibitem{Mathai1}
Varghese~Mathai and Danny~Stevenson,
Chern character in twisted $K$-theory: Equivariant and holomorphic cases,
{\it Commun. Math. Phys.} {\bf 236}:161-186, 2003.
%
\bibitem{McDuff1}
Dusa~McDuff,
The moment map for circle actions on symplectic manifolds,
{\it J. Geom. Phys.} {\bf 5}(2):149--160, 1988.
%
\bibitem{Meinrenken}
Eckhard~Meinrenken,
Symplectic surgery and the Spin-c Dirac operator,
{\it Advances in Mathematics}, {\bf 134}  (1998), 240--277.
\bibitem{Noda}
Tomonori~Noda,
Reduction of locally conformal symplectic manifolds with examples of non-K\"ahler manifolds,
{\it Tsukuba J. Math.}, {\bf 28}(1):127-136, 2004.
%
\bibitem{Nunes}
Joana~M.~Nunes da Costa,
Reduction of complex Poisson manifolds,
{\it Portugal. Math.}, {\bf 54}(4):467--476, 1997.
%
\bibitem{Ornea}
Liviu~Ornea and Paolo~Piccinni,
Cayley $4$-frames and a quaternion K\"ahler reduction related to ${\rm Spin}(7)$,
In {\it Global differential geometry: the mathematical legacy of Alfred Gray (Bilbao, 2000)}, 401--405,
Contemp. Math., 288,
Amer. Math. Soc., Providence, RI, 2001. 
%
\bibitem{Ortega}
Juan-Pablo~Ortega and Tudor~S.~Ratiu,
Singular reduction of Poisson manifolds,
{\it Lett. Math. Phys.} {\bf 46}(2):359--372, 1998.
%
\bibitem{Ortega1}
Juan-Pablo~Ortega and Tudor~S.~Ratiu,
Momentum maps and Hamiltonian reduction,
{\it Progess in Mathematics} {\bf 222}, Birkh\"auser-Verlag, Boston, 2004.
%
\bibitem{Proudfoot}
Nicholas~J.~Proudfoot,
Hyperkahler analogues of Kahler quotients,
Ph. D. thesis, U.C. Berkeley, Spring 2004,
{\it math.AG/0405233}.
%
\bibitem{Severa}
Pavol~\v Severa,
Letters to Weinstein, 
{\it http://sophia.dtp.fmph.uniba.sk/$\sim$severa/letters/no1.ps}, 1998.
%
\bibitem{Stienon}
Mathieu~Sti\'enon and Ping~Xu,
Reduction of generalized complex structures,
{\it math.DG/0509393}, 2005.
%
\bibitem{Uchino}
Kyousuke~Uchino,
Remarks on the definition of a Courant algebroid.
{\it Lett. Math. Phys.,} {\bf 60}(2):171-175, 2002.
%
\bibitem{Vaisman}
Izu~Vaisman,
Reduction and submanifolds of generalized complex manifolds.
{\it math.DG/0511013}, 2005.
%
\bibitem{Willett}
Christopher~Willett,
Contact reduction,
{\it Trans. Amer. Math. Sco.}, {\bf 354}(10):4245-4260, 2002.
\end{thebibliography}
\end{document}